\documentclass[11pt]{article}
\usepackage{amscd}
\usepackage{amsfonts}
\usepackage{amsmath}
\usepackage{amssymb}
\usepackage{amsthm}
\usepackage{bbm}
\usepackage{CJK}
\usepackage{fancyhdr}
\usepackage{graphicx}
\usepackage{hyperref}
\usepackage{indentfirst}
\usepackage{latexsym}
\usepackage{mathrsfs}
\usepackage{xypic}
\allowdisplaybreaks

\usepackage[top=1in,bottom=1in,left=1.25in,right=1.25in]{geometry}
\def\<{\langle}
\def\>{\rangle}
\def\a{\alpha}
\def\b{\beta}

\def\c{\cdot}

\def\D{\Delta}

\def\g{\gamma}

\def\lr{\longrightarrow}

\def\o{\otimes}

\def\v{\epsilon }

\date{}
\begin{document}
\renewcommand{\baselinestretch}{1.2}
\renewcommand{\arraystretch}{1.0}
\title{\bf The Hom-Long dimodule category\\ and nonlinear equations}
\author{{\bf Shengxiang Wang$^{1}$, Xiaohui Zhang$^{2}$,
        Shuangjian Guo$^{3}$\footnote
        {Correspondence: shuangjianguo@126.com} }\\
1.~ School of Mathematics and Finance, Chuzhou University,\\
 Chuzhou 239000,  China \\
2.~  School of Mathematical Sciences, Qufu Normal University, \\Qufu Shandong 273165, China.\\
 3.~ School of Mathematics and Statistics, Guizhou University of\\ Finance and Economics, Guiyang 550025, China.}
 \maketitle
\begin{center}
\begin{minipage}{13.cm}

{\bf \begin{center} ABSTRACT \end{center}}
In this paper, we construct a kind of new braided monoidal category over two Hom-Hopf algerbas  $(H,\alpha)$  and $(B,\beta)$
 and associate it with two nonlinear equations.
We first introduce the notion of an $(H,B)$-Hom-Long dimodule
and show that the Hom-Long dimodule category  $^{B}_{H} \Bbb L$   is an autonomous category.
Second, we prove that the category  $^{B}_{H} \Bbb L$  is a braided monoidal category if $(H,\alpha)$ is  quasitriangular
and $(B,\beta)$ is coquasitriangular and get a solution of the quantum Yang-Baxter equation.
Also, we  show that the category  $^{B}_{H} \Bbb L$ can be viewed as  a subcategory of the Hom-Yetter-Drinfeld category $^{H\o B}_{H\o B} \Bbb {HYD}$.
Finally, we obtain a solution of the Hom-Long equation from the Hom-Long dimodules.\\

{\bf Key words}:   Hom-Long dimodule; Hom-Yetter-Drinfeld category; Yang-Baxter equation; Hom-Long equation.\\

 {\bf 2010 Mathematics Subject Classification:} 16A10; 16W30
 \end{minipage}
 \end{center}
 \normalsize\vskip1cm

\section*{INTRODUCTION}
\def\theequation{0. \arabic{equation}}
\setcounter{equation} {0}

The study of Hom-algebras can be traced back to  Hartwig,  Larsson and  Silvestrov's work in \cite{Hartwig},
 where the notion of Hom-Lie algebra in the context of q-deformation theory of Witt and Virasoro algebras \cite{Hu} was introduced,
which plays an important role in physics, mainly in conformal field theory.
Hom-algebras and Hom-coalgebras were introduced by Makhlouf and Silvestrov  \cite{Makhlouf2008} as
generalizations of ordinary algebras and coalgebras in the following sense:
the associativity of the multiplication is replaced by
the Hom-associativity and similar for Hom-coassociativity.
They also defined the structures of Hom-bialgebras and Hom-Hopf algebras,
and described some of their properties extending properties of ordinary bialgebras and Hopf algebras
in \cite{Makhlouf2009, Makhlouf2010}.
In \cite{Caenepeel2011}, Caenepeel and Goyvaerts studied Hom-bialgebras and Hom-Hopf algebras from a categorical view point,
and called them monoidal Hom-bialgebras and monoidal Hom-Hopf algebras respectively,
 which are different from the normal Hom-bialgebras and Hom-Hopf algebras in \cite{Makhlouf2009}.
Many more properties and structures of Hom-Hopf algebras have been developed,
see  \cite{CX14, Gohr2010, ma2014, ZD} and references cited therein.

Later, Yau \cite{Yau2009, Yau3} proposed the definition of quasitriangular Hom-Hopf algebras
and showed that each quasitriangular Hom-Hopf algebra yields a solution of the Hom-Yang-Baxter equation.
The  Hom-Yang-Baxter equation reduces to the usual Yang-Baxter equation when the twist map is trivial.
Several classes of solutions of the  Hom-Yang-Baxter equation were constructed from different respects,
 including those associated
to Hom-Lie algebras \cite{Fang, wang&GUO2020, Yau2009, Yau2011}, Drinfelds (co)doubles \cite{CWZ2014, ZGW, ZWZ},  and Hom-Yetter-Drinfeld modules \cite{chen2014, LIU2014, ma2017, ma2019, Makhlouf2014, WCZ2012, YW}.

It is well-known that classical nonlinear equations in Hopf algebra theory including the quantum Yang-Baxter equation,
 the Hopf equation, the pentagon equation, and the Long equation.
In \cite{Militaru}, Militaru proved that each Long dimodule gave rise to a solution for the Long equation.
Long dimodules are the building stones of the Brauer-Long group.
In the case where $H$ is commutative, cocommutative and faithfully projective,
 the Yetter-Drinfeld category $^{H}_{H}\Bbb{YD}$ is precisely the Long dimodule category $^{H}_{H}\Bbb{L}$.
Of course, for an arbitrary $H$, the categories $^{H}_{H}\Bbb{YD}$ and $^{H}_{H}\Bbb{L}$  are basically different.
In \cite{CWZ2013}, Chen  et al. introduced the concept of Long dimodules over a monoidal Hom-bialgebra and discussed its relation with Hom-Long equations.
Later, we \cite{wang&ding2017} extended Chen's work to generalized Hom-Long dimodules over monoidal Hom-Hopf algebras and obtained a kind solution for the quantum Yang-Baxter equation.
For more details about  Long dimodules, see \cite{Long, Lu, WangShuanhong, Zhang}  and references cited therein.

The main purpose of this paper is to construct a new braided monoidal category and present solutions for two kinds of nonlinear equations.
Different to our previous work in \cite{wang&ding2017}, in the present paper we do all the work over  Hom-Hopf algebras,
which is more unpredictable than the monoidal version.
Since Hom-Hopf algebras and monoidal Hom-Hopf algebras are different concepts, it turns out that our definitions,
formulas and results are also different from the ones in \cite{wang&ding2017}.
Most important, we  associate quantum Yang-Baxter equations and Hom-Long equations to the  Hom-Long dimodule categories.

This paper is organized as follows.
In Section 1, we recall some basic definitions about Hom-(co)modules and (co)quasitriangular Hom-Hopf algebras .

In Section 2, we first introduce the notion of $(H,B)$-Hom-Long dimodules over Hom-bialgebras $(H,\alpha)$ and $(B, \b)$,
then we show that the   Hom-Long dimodule category $^{B}_{H}\Bbb{L}$ forms an autonomous category (see Theorem 2.6)
and prove that the category is equivalent to the category of left $B^{\ast op} \o H$-Hom-modules (see Theorem 2.7).

In Section 3, for a quasitriangular Hom-Hopf algebra $(H,R, \alpha)$ and a coquasitriangular Hom-Hopf algebra  $(B,\langle|\rangle,\beta)$,
we prove that the Hom-Long dimodule category $^{B}_{H}\Bbb{L}$ is a subcategory of the
Hom-Yetter-Drinfeld category $^{H\o B}_{H\o B}\Bbb{HYD}$ (see Theorem 3.5),
and show that  the braiding yields a solution for the quantum Yang-Baxter equation (see Corollary  3.2).

In Section 4, we prove that the category $_{H}\Bbb{M}$ over a triangular Hom-Hopf algebra
(resp., $^{H}\Bbb{M}$ over a cotriangular Hom-Hopf algebra) is a
Hom-Long dimodule subcategory of  $^{B}_{H}\Bbb{L}$ (see Propositions 4.1 and 4.2).
We also show that the  Hom-Long dimodule category  $^{B}_{H}\Bbb{L}$
is symmetric in case $(H,R, \alpha)$ is triangular and $(B,\langle|\rangle,\beta)$ is cotriangular  (see Theorem 4.3).

In Section 5, we introduce the notion of  $(H,\a)$-Hom-Long dimodules and  obtain a solution for the Hom-Long equation (see Theorem 5.10).

\section{PRELIMINARIES}
\def\theequation{\arabic{section}.\arabic{equation}}
\setcounter{equation} {0}

Throughout this paper, $k$ is a fixed field.
 Unless otherwise stated, all vector spaces, algebras, modules, maps and unadorned tensor products are over $k$.
 For a coalgebra $C$, the coproduct  will be denoted by $\Delta$.
 We adopt a Sweedler's notation $\triangle(c)=c_{1}\otimes c_{2}$, for any $c\in C$, where the summation is understood.
We refer to
\cite{R, Sweedler}
 for the Hopf algebra theory and terminology.

 We now recall some useful definitions in \cite{Li2014, Makhlouf2008, Makhlouf2009,  Makhlouf2010, Yau1, Yau3}.
\smallskip

 {\bf Definition 1.1.} A Hom-algebra is a quadruple $(A,\mu,1_A,\a)$ (abbr. $(A,\a)$), where $A$ is a $k$-linear space,
  $\mu: A\o A \lr A$ is a $k$-linear map, $1_A \in A$ and $\a$ is an automorphism of $A$, such that
 \begin{eqnarray*}
 &(HA1)& \a(aa')=\a(a)\a(a');~~\a(1_A)=1_A,\\
 &(HA2)& \a(a)(a'a'')=(aa')\a(a'');~~a1_A=1_Aa=\a(a)
 \end{eqnarray*}
 are satisfied for $a, a', a''\in A$. Here we use the notation $\mu(a\o a')=aa'$.
\smallskip

 {\bf Definition 1.2.}
 Let $(A,\alpha)$ be a Hom-algebra.
 A left $(A,\alpha)$-Hom-module is a triple $(M,\rhd,\nu)$, where $M$ is a linear space,
 $\rhd: A\o M \lr M$ is a linear map, and $\nu$ is an automorphism of $M$, such that
  \begin{eqnarray*}
 &(HM1)& \nu(a\rhd m)=\alpha(a)\rhd \nu(m),\\
 &(HM2)& \alpha(a)\rhd (a'\rhd m)=(aa')\rhd \nu(m);~~ 1_A\rhd m=\nu(m)
 \end{eqnarray*}
 are satisfied for $a, a' \in A$ and $m\in M$.

 Let $(M,\rhd_M,\nu_M)$ and $(N,\rhd_N,\nu_N)$ be two left $(A,\alpha)$-Hom-modules.
 Then a linear morphism $f: M\lr N$ is called a morphism of left $(A,\alpha)$-Hom-modules if
 $f(h\rhd_M m)=h\rhd_N f(m)$ and $\nu_N\circ f=f\circ \nu_M$.
 \smallskip

 {\bf Definition 1.3.}  A Hom-coalgebra is a quadruple $(C,\D,\v,\b)$ (abbr. $(C,\b)$),
  where $C$ is a $k$-linear space, $\D: C \lr C\o C$, $\v: C\lr k$ are $k$-linear maps,
  and $\b$ is an automorphism of $C$, such that
 \begin{eqnarray*}
 &(HC1)& \b(c)_1\o \b(c)_2=\b(c_1)\o \b(c_2);~~\v\circ \b=\v;\\
 &(HC2)& \b(c_{1})\o c_{21}\o c_{22}=c_{11}\o c_{12}\o \b(c_{2});~~\v(c_1)c_2=c_1\v(c_2)=\b(c)
 \end{eqnarray*}
 are satisfied for $c \in C$.
\smallskip

 {\bf Definition 1.4.}
  Let $(C,\b)$ be a Hom-coalgebra.
  A left $(C,\b)$-Hom-comodule is a triple $(M,\rho,\mu)$, where $M$ is a linear space,
   $\rho:  M\lr C\o M$ (write $\rho(m)=m_{(-1)}\o m_{(0)},~\forall m\in M$) is a linear map,
 and $\mu$ is an automorphism of $M$, such that
 \begin{eqnarray*}
 &(HCM1)&\mu(m)_{(-1)}\o \mu(m)_{(0)}=\b(m_{(-1)})\o \mu(m_{(0)}),~\v(m_{(-1)})m_{(0)}=\mu(m);\\
 &(HCM2)&\b(m_{(-1)})\o m_{(0)(-1)}\o m_{(0)(0)}= m_{(-1)1}\o m_{(-1)2}\o \mu(m_{(0)})
 \end{eqnarray*}
 are satisfied for all $m\in M$.

  Let $(M,\rho^M,\mu_M)$ and $(N,\rho^N,\mu_N)$ be two left $(C,\b)$-Hom-comodules.
   Then a linear map $f: M\lr N$ is called a map of left $(C,\b)$-Hom-comodules if
   $f(m)_{(-1)}\o f(m)_{(0)}=m_{(-1)}\o f(m_{(0)})$ and $\mu_N\circ f=f\circ \mu_M$.

 \smallskip

 {\bf Definition 1.5.}  A Hom-bialgebra is a sextuple $(H,\mu,1_H,\D,\v,\g)$ (abbr. $(H,\g)$), where
 $(H,\mu,1_H,\g)$ is a Hom-algebra and $(H,\D,\v,\g)$ is a Hom-coalgebra, such that
 $\D$ and $\v$ are morphisms of Hom-algebras, i.e.
 $$
 \D(hh')=\D(h)\D(h');~\D(1_H)=1_H\o 1_H;~
 \v(hh')=\v(h)\v(h');~\v(1_H)=1.
 $$
 Furthermore, if there exists a linear map $S: H\lr H$ such that
 $$
 S(h_1)h_2=h_1S(h_2)=\v(h)1_H~\hbox{and}~S(\g(h))=\g(S(h)),
 $$
 then we call $(H,\mu,1_H,\D,\v,\g,S)$ (abbr. $(H,\g,S)$) a Hom-Hopf algebra.
 \smallskip

 {\bf Definition 1.6.} (\cite{Li2014}) Let $(H, \b)$ be a Hom-bialgebra,
 $(M,\rhd, \mu)$ a  left $(H,\b)$-module with action $\rhd: H\o M\lr M, h\o m\mapsto h\rhd m$
 and $(M,\rho, \mu)$ a left $(H,\b)$-comodule with coaction $\rho: M\lr H\o M, m\mapsto m_{(-1)}\o m_{(0)}$.
  Then we call $(M,\rhd, \rho, \mu)$ a (left-left) Hom-Yetter-Drinfeld module over $(H,\b)$ if the following condition holds:
 $$
 (HYD)~~~h_1\b(m_{(-1)})\o (\b^3(h_2)\rhd m_{(0)}=(\b^2(h_1)\rhd m)_{(-1)}h_2\o (\b^2(h_1)\rhd m)_{(0)},
 $$
 where $h\in H$ and $m\in M$.

 When $H$ is a Hom-Hopf algebra, then the condition $(HYD)$ is equivalent to
 $$
 (HYD)' ~~\rho(\b^4(h)\rhd m)=\b^{-2}(h_{11}\b(m_{(-1)}))S(h_2)\o (\b^3(h_{12})\rhd m_0).
 $$

{\bf Definition 1.7.} (\cite{Li2014}) Let $(H, \b)$ be a Hom-bialgebra.
  A Hom-Yetter-Drinfeld category $^{H}_{H}\Bbb Y\Bbb D$
 is a pre-braided monoidal category whose objects are left-left Hom-Yetter-Drinfeld modules, morphisms are both
 left $(H, \beta)$-linear and $(H, \beta)$-colinear maps, and its pre-braiding $C_{-, -}$  is given by
\begin{eqnarray}
 C_{M, N} (m \o n) = \beta^{2}(m_{(-1)})\rhd\nu^{-1}(n)\o \mu^{-1}(m_{0}),
\end{eqnarray}
for all $m \in (M,\mu)\in{}^{H}_{H}\Bbb Y\Bbb D$ and $n \in (N,\nu) \in{}^{H}_{H}\Bbb Y\Bbb D$.
 \smallskip

{\bf Definition 1.8.}
A quasitriangular Hom-Hopf algebra is a octuple $(H, \mu, 1_{H},\Delta, \epsilon, S, \beta, R)$ (abbr. $(H, \beta, R)$)
in which $(H, \mu, 1_{H}, \Delta, \epsilon, S, \beta)$ is a Hom-Hopf algebra and
$R=R^{(1)}\o R^{(2)}\in H\o H$,  satisfying the following axioms (for all $h\in H$ and $R=r$):
\begin{eqnarray*}
&&(QHA1)~\epsilon(R^{(1)})R^{(2)}=R^{(1)}\epsilon(R^{(2)})=1_{H},\\
&&(QHA2)~\Delta(R^{(1)})\otimes\beta(R^{(2)})=\beta(R^{(1)})\otimes\beta(r^{(1)})\otimes R^{(2)}r^{(2)},\\
&&(QHA3)~\beta(R^{(1)})\otimes\Delta(R^{(2)})=R^{(1)}r^{(1)}\otimes \beta(r^{(2)})\otimes \beta(R^{(2)}),\\
&&(QHA4)~\Delta^{cop}(h)R=R\Delta(h),\\
&&(QHA5)~\beta(R^{(1)})\o\beta(R^{(2)})=R^{(1)}\o R^{(2)},
\end{eqnarray*}
where $\Delta^{cop}(h)=h_{2}\otimes h_{1}$ for all $h\in H$.
A quasitriangular Hom-Hopf algebra $(H,R,\beta)$ is called triangular if
$R^{-1}=R^{(2)}\otimes R^{(1)}.$
\medskip

{\bf Definition 1.9.}
A coquasitriangular Hom-Hopf algebra  is a Hom-Hopf algebra $(H, \beta)$
together with a bilinear form $\langle|\rangle$ on $(H, \beta)$ (i.e. $\langle|\rangle\in$ Hom($H\otimes H,k$))
such that the following axioms hold:
\begin{eqnarray*}
&&(CHA1)~\langle hg|\beta(l)\rangle=\langle\beta(h)|l_{2}\rangle\langle \beta(g)|l_{1}\rangle,\\
&&(CHA2)~\langle \beta(h)|gl\rangle=\langle h_{1}|\beta(g)\rangle\langle h_{2}|\beta(l)\rangle,\\
&&(CHA3)~\langle h_{1}|g_{1}\rangle g_{2}h_{2}=h_{1}g_{1}\langle h_{2}|g_{2}\rangle,\\
&&(CHA4)~\langle 1|h\rangle=\langle h|1\rangle=\epsilon(h),\\
&&(CHA5)~\langle \beta(h)|\beta(g)\rangle=\langle h|g\rangle
\end{eqnarray*}
for all $h,g,l\in H$.
A coquasitriangular Hom-Hopf algebra $(H,\langle|\rangle, \beta)$ is called cotriangular if $\langle|\rangle$
 is convolution invertible in the sense of
$
\langle h_{1}|g_{1}\rangle\langle g_{2}|h_{2}\rangle=\epsilon(h)\epsilon(g),
$
for all $h,g\in H$.

\section{  Hom-Long dimodules over Hom-bialgebras}
\def\theequation{\arabic{section}.\arabic{equation}}
\setcounter{equation} {0}

In this section, we will introduce the notion of  Hom-Long dimodules
and prove that the Hom-Long dimodule  category is an autonomous category.
\medskip

\noindent{\bf Definition 2.1.}
Let  $(H,\alpha)$ and  $(B,\beta)$ be two  Hom-bialgebras.
A left-left $(H,B)$-Hom-Long dimodule is a quadrupl  $(M,\c, \rho,\mu)$,
where $(M, \c, \mu)$ is a left  $(H,\alpha)$-Hom-module
and  $(M, \rho, \mu)$ is a left  $(B,\beta)$-Hom-comodule such that
\begin{eqnarray}
\rho(h\c m)=\beta(m_{(-1)})\o \alpha(h)\c m_{(0)},
\end{eqnarray}
for all $h\in H$ and $m\in M$.
We denote by  $^{B}_{H} \Bbb L$ the category of left-left $(H,B)$-Hom-Long dimodules,
morphisms being $H$-linear $B$-colinear maps.
\medskip

\noindent{\bf Example 2.2.}
Let  $(H,\alpha)$ and  $(B,\beta)$ be two   Hom-bialgebras.
Then $(H\o B,\alpha\o\beta)$ is an $(H,B)$-Hom-Long dimodule
with left $(H,\alpha)$-action $h\c(g\o x)=hg\o \beta(x)$ and  left $(B,\beta)$-coaction
$\rho(g\o x)=x_{1}\o(\alpha(g)\o x_{2})$, where $h,g\in H, x\in B$.
\medskip

\noindent{\bf Proposition 2.3.}
Let $(M, \mu), (N, \nu)$ be two $(H,B)$-Hom-Long dimodules,
then $(M\o N, \mu\o\nu)$  is an  $(H,B)$-Hom-Long dimodule with structures:
\begin{eqnarray*}
&&h\c (m\o n)=h_{1}\c m\o h_{2}\c n,\\
&&\rho(m\o n)=\b^{-2}(m_{(-1)}n_{(-1)})\o m_{(0)}\o n_{(0)},
\end{eqnarray*}
for all $m\in M, n\in N$ and $h\in H$.
\medskip

\noindent{\bf Proof.}
From Theorem 4.8 in \cite{ma2017},
$(M\o N, \mu\o\nu)$  is both a left  $(H,\alpha)$-Hom-module and a left  $(B,\beta)$-Hom-comodule.
It remains to check that the compatibility condition (2.1) holds.
For any $m\in M, n\in N$ and $h\in H$, we have
\begin{eqnarray*}
\rho(h\c (m\o n))
&=&\b((h_{1}\c m)_{(-1)}(h_{2}\c n)_{(-1)})\o(h_{1}\c m)_{(0)}\o(h_{2}\c n)_{(0)}\\
&=&\beta^{-1}(m_{(-1)}n_{(-1)})\o\alpha(h_{1})\c m_{(0)}\o\alpha(h_{2})\c n_{(0)}\\
&=&\beta((m\o n)_{(-1)})\o\alpha(h)\c ((m\o n)_{(0)}),
\end{eqnarray*}
as desired. This completes the proof. \hfill $\square$
\medskip

\noindent{\bf Proposition 2.4.}
The Hom-Long dimodule category $^{B}_{H} \Bbb L$
is a monoidal category, where the tensor product is given in Proposition 2.3,
the unit $I=(k,id)$, the associator and the constraints are given as follows:
\begin{eqnarray*}
&&a_{U,V,W}: (U\o V)\o W\rightarrow U\o (V\o W), (u\o v)\o w\rightarrow \mu^{-1}(u)\o(v\o \omega(w)),\\
&&l_{V}: k\o V\rightarrow V, k\o v\rightarrow k\nu(v), r_{V}: V\o k\rightarrow V, v\o k \rightarrow k\nu(v),
\end{eqnarray*}
for $u\in (U,\mu)\in {}^{B}_{H} \Bbb L, v\in (V,\nu)\in{}^{B}_{H} \Bbb L, w\in (W,\omega)\in{} ^{B}_{H} \Bbb L.$
\medskip

\noindent{\bf Proof.} Straightforward. \hfill $\square$

\medskip

\noindent{\bf Proposition 2.5.}
Let $H$ and $B$ be two Hom-Hopf algebras with bijective antipodes.
For any   Hom-Long dimodule $(M,\mu)$ in $^{B}_{H} \Bbb L$,
set $M^\ast = Hom_k(M,k)$, with the $(H,\alpha)$-Hom-module and the $(B,\b)$-Hom-comodule structures:
\begin{eqnarray*}
&&\theta_{M^\ast}:  H\o M^\ast\longrightarrow M^\ast,~~
 (h\cdot f )(m)=f(S_H\a^{-1}(h)\cdot\mu^{-2}(m)),\\
&&\rho_{M^\ast}:  M^\ast \longrightarrow B \o M^\ast, ~~
 f_{(-1)} \o f_{(0)}(m)=S^{-1}_B\b^{-1}(m_{(-1)}) \o f(\mu^{-2}(m_{(0)})),
\end{eqnarray*}
and the Hom-structure map $\mu^{\ast}$ of $M^\ast$ is $\mu^{\ast}(f)(m)=f(\mu^{-1}(m))$.
Then $M^\ast$ is an object in $^{B}_{H} \Bbb L$. Moreover, $^{B}_{H} \Bbb L$ is a left autonomous category.
\medskip

\noindent{\bf Proof.}
It is not hard to check that $(M^\ast,\theta_{M^\ast},\mu^{\ast})$ is an $(H,\alpha)$-Hom-module
 and $(M^\ast,\rho_{M^\ast},\mu^\ast)$ is a $(B,\beta)$-Hom-comodule.
Further, for any $f \in M^\ast$, $m \in M$, $h\in H$, we have
\begin{eqnarray*}
(h\c f)_{(-1)} \o (h\c f)_{(0)}(m)
&=&S^{-1}_B\b^{-1}(m_{(-1)}) \o (h\c f)(\mu^{-2}(m_{(0)}))\\
&=&S^{-1}_B\b^{-1}(m_{(-1)}) \o f(S_H\a^{-1}(h) \c \mu^{-4}(m_{(0)})),\\
\b(f_{(-1)}) \o (\alpha(h)\c f_{(0)})(m)
&=&\b(f_{(-1)}) \o f_{(0)}(S_H(h) \c \mu^{-2}(m) )\\
&=&\b(S^{-1}_B\b^{-2}(m_{(-1)} )) \o f( \mu^{-2}( S_H \alpha (h) \c \mu^{-2}(m_{(0)}) )) \\
&=&S^{-1}_B\b^{-1}(m_{(-1)}) \o f(S_H\a^{-1}(h) \c \mu^{-4}(m_{(0)})).
\end{eqnarray*}
Thus $M^\ast \in {}^{B}_{H} \Bbb L$.

Moreover, for any $f \in M^\ast$ and $m\in M$, one can define the left evaluation map and the left coevaluation map by
\begin{eqnarray*}
ev_M:   f \otimes  m \longmapsto  f(m),
~coev_M:  1_k \longmapsto  \sum e_i \o e^i,
\end{eqnarray*}
where $e_i$ and $e^i$ are dual bases in $M$ and $M^\ast$ respectively.
Next, we will show $(M^\ast,ev_M,coev_M)$ is the left dual of $M$.

It is easy to see that $ev_M$ and $coev_M$ are morphisms in $^{B}_{H} \Bbb L$.
For this, we need the following computation
\begin{eqnarray*}
&&~~~(r_M\circ(id_M \o ev_M)\circ a_{M,M^\ast,M} \circ (coev_M \o id_M)\circ l_M^{-1})(m)\\
&&=(r_M\circ(id_M \o ev_M)\circ a_{M,M^\ast,M})(\sum_{i} (e_i \o   e^i) \o \mu^{-1}(m))\\
&&=(r_M\circ(id_M \o ev_M))(\sum_{i} \mu^{-1}(e_i) \o  ( e^i \o m))\\
&&=r_M(\sum_{i} \mu^{-1}(e_i) \o e^i(m ) )\\
&&= r_M(\mu^{-1}(m) \o 1_k )=m.
\end{eqnarray*}
Similarly, we get
\begin{eqnarray*}
&&~~~(l_{M^\ast}\circ (ev_M \o id_{M^\ast})\circ a^{-1}_{M^\ast,M,M^\ast}\circ (id_{M^\ast} \o coev_M)\circ r_{M^\ast}^{-1})(f)\\
&&=(l_{M^\ast}\circ (ev_M \o id_{M^\ast})\circ a^{-1}_{M^\ast,M,M^\ast})(\sum_{i} {\mu}^{\ast-1}(f) \o (e_i \o e^i)) \\
&&=(l_{M^\ast}\circ (ev_M \o id_{M^\ast}))(\sum_{i} f\o e_i) \o \mu^{\ast-1}(e^i) )\\
&&=l_{M^\ast}(\sum_{i}  f(e_i)   \o \mu^{\ast-1}(e^i) )\\
&&=l_{M^\ast}(1_k \o {\mu}^{\ast-1}(f) )=f.
\end{eqnarray*}
So $^{B}_{H} \Bbb L$ admits the left duality. The proof is finished.\hfill $\square$
\medskip

\noindent{\bf Theorem 2.6.}
The Hom-Long dimodule category $^{B}_{H} \Bbb L$ is an autonomous category.
\medskip

\noindent{\bf Proof}
By Proposition 2.5, it is sufficient to show that $^{B}_{H} \Bbb L$ is also a right autonomous category.
In fact, for any $(M,\mu) \in {}^{B}_{H} \Bbb L$, its right dual $({}^\ast M, \widetilde{coev}_M, \widetilde{ev}_M)$ is defined as follows:

$\bullet$ ${}^\ast M = Hom_k(M,k)$ as $k$-modules, with the Hom-module and Hom-comodule structures:
\begin{eqnarray*}
&( h\cdot f)(m)=f(S^{-1}_H\a^{-1}(h)\cdot\mu^{-2}(m)),\\
&f_{(-1)} \o f_{(0)}(m)=S_B\b^{-1}(m_{(-1)}) \o f(\mu^{-2} (m_{(0)})),
\end{eqnarray*}
where $f \in {}^\ast M$, $m \in M$,
and the Hom-structure map $\mu^{\ast}$ of $^\ast M$ is $\mu^{\ast}(f)(m)=f(\mu^{-1}(m))$;

$\bullet$ The right evaluation map and the right coevaluation map are given by
\begin{eqnarray*}
\widetilde{ev}_M:   m\otimes f\longmapsto  f(m),
~\widetilde{coev}_M:  1_k \longmapsto  \sum a^i \o a_i.
\end{eqnarray*}
where $a_i$ and $a^i$ are dual bases of $M$ and ${}^\ast M$ respectively.
By similar verification in Proposition 2.5, one may check that
$^{B}_{H} \Bbb L$ is a right autonomous category, as required.
This completes the proof.\hfill $\square$
\medskip

Recall from \cite{ZGW} that for any finite dimensional Hom-Hopf algebra $B$, $B^\ast$ is also a Hom-Hopf algebra with the following structures
\begin{eqnarray*}
&(f\ast g)(y):=f(\b^{-2}(y_1))g(\b^{-2}(y_2)),~~\Delta_{B^\ast}(f)(xy):=f(\b^{-2}(xy)),\\
&1_{B^\ast}:=\v,~~\v_{B^\ast}(f):=f(1_H),~~S_{B^\ast}:=S^\ast,~~\a_{B^\ast}(f):=f\circ \b^{-1},
\end{eqnarray*}
where $x,y \in H$, $f,g \in B^{\ast}$.

\noindent{\bf Theorem 2.7.}
If $B$ is a finite dimensional Hom-Hopf algebra, then the  Hom-Long dimodule category $^{B}_{H} \Bbb L$ is identified to the category of left $B^{\ast op} \o H$-Hom-modules, where $B^{\ast op} \o H$ means the usual tensor product Hom-Hopf algebra.
\medskip

\noindent{\bf Proof}
Define the functor $\Psi$ from ${}_{B^{\ast op} \o H} \Bbb{M}$ to $^{B}_{H} \Bbb L$ by
\begin{eqnarray*}
\Psi(M):=M\mbox{~as~}k\mbox{-module~},~~~~\Psi(f):=f,
\end{eqnarray*}
where $(M, \mu, \rightharpoondown)$ is a $B^{\ast op} \o H$-Hom-module, $f:M\rightarrow N$ is a morphism of $B^{\ast op} \o H$-Hom-modules.
Further, the $H$-action on $M$ is defined by
\begin{eqnarray*}
h \c m:= (\v_B \o h)\rightharpoondown m, ~~~~\mbox{for~all~}m\in M,~~h\in H,
\end{eqnarray*}
and the $B$-coaction on $M$ is given by
\begin{eqnarray*}
m_{(-1)} \o m_{(0)}:= \sum e_i \o (e^i \o 1_H)\rightharpoondown m,
\end{eqnarray*}
where $e_i$ and $e^i$ are dual bases of $B$ and $B^\ast$ respectively.

First, we will show $(M,\mu,\c)$ is a left $(H,\a)$-Hom-module. Actually, for any $m \in M$, $h,g \in H$, we have
$
1_H \c m = (\v_B \o 1_H)\rightharpoondown m = \mu(m),
$
and
\begin{eqnarray*}
\a(h)\c(g \c m) &=& (\v_B \o \a(h)) \rightharpoondown ( (\v_B \o g) \rightharpoondown m ) \\
&=& (\v_B \o hg)\rightharpoondown \mu(m)=(hg) \c \mu(m),
\end{eqnarray*}
which implies $(M,\mu,\c) \in {}_H \Bbb{M}$.

Second, one can show that $(M,\mu) \in {}^B \Bbb{M}$ in a similar way.

At last, for any $m \in M$, $h \in H$, we have
\begin{eqnarray*}
 (h \c m)_{(-1)} \o (h \c m)_{(0)}
&=& \sum e_i \o (e^i \o 1_H) \rightharpoondown (h \c m) \\
&=& \sum e_i \o ( e^i \o \a(h) ) \rightharpoondown \mu(m) \\
&=& \sum \b(e_i) \o ( (\v_B \o 1_H) (e^i \o h ) \rightharpoondown \mu(m) \\
&=& \sum \b(e_i) \o ( (\v_B \o h ) (e^i \o 1_H ) \rightharpoondown \mu(m) \\
&=& \sum \b(e_i) \o \a(h) \c ((e^i \o 1_H ) \rightharpoondown \mu(m))\\
&=& \b(m_{(-1)}) \o \a(h) \c m_{(0)},
\end{eqnarray*}
which implies $(M,\mu) \in {}^B_H\Bbb{L}$.

Conversely, for any object $(M,\mu)$, $(N,\nu)$, and morphism $f:U\rightarrow V$ in ${}^B_H\Bbb{L}$, one can define a functor $\Phi$ from ${}^B_H\Bbb{L}$ to ${}_{B^{\ast op} \o H} \Bbb{M}$
\begin{eqnarray*}
\Phi (M):=M\mbox{~as~}k\mbox{-modules~},~~~~\Phi(f):=f,
\end{eqnarray*}
where the $(B^{\ast op}\o H,\b^*\o \a)$-Hom-module structure on $M$ is given by
\begin{eqnarray*}
(p \o h)\rightharpoondown m = p(m_{(-1)}) h \c \mu^{-1}(m_{(0)}),
\end{eqnarray*}
for all $p \in B^{\ast},~h \o H,~m \in M.$
It is straightforward to check that $(M,\mu, \rightharpoondown)$ is an object in $^{B}_{H} \Bbb L$ to ${}_{B^{\ast op} \o H} \Bbb{M}$, and hence $\Phi$ is well defined.

Note that $\Phi$ and $\Psi$ are inverse with each other. Hence the conclusion holds.

\section{New braided momoidal categories over Hom-Long dimodules}
\def\theequation{\arabic{section}.\arabic{equation}}
\setcounter{equation} {0}

In this section, we will prove that the Hom-Long dimodule category $^{B}_{H} \Bbb L$
over a quasitriangular Hom-Hopf algebra $(H,R,\alpha)$
and a coquasitriangular Hom-Hopf algebra $(B,\langle|\rangle,\beta)$
is a braided  monoidal subcategory
 of the Hom-Yetter-Drinfeld category $^{H\o B}_{H\o B} \Bbb{HYD}$.
 \medskip

\noindent{\bf Theorem 3.1.}
Let $(H,R,\alpha)$ be a quasitriangular Hom-Hopf algebra
and $(B,\langle|\rangle,\beta)$ a coquasitriangular Hom-Hopf algebra.
Then the category $^{B}_{H} \Bbb L$  is a braided  monoidal category with braiding
\begin{eqnarray}
C_{M,N}: M\o N\rightarrow N\o M, m\o n\rightarrow \langle m_{(-1)}|n_{(-1)}\rangle
                                      R^{(2)}\c \nu^{-2}(n_{(0)})\o R^{(1)}\c \mu^{-2}(m_{(0)}),
\end{eqnarray}
for all $m\in (M,\mu)\in{}^{B}_{H} \Bbb L$ and $ n\in (N,\nu)\in{}^{B}_{H} \Bbb L.$

\medskip

\noindent{\bf Proof.}
We will first show that the braiding  $C_{M,N}$ is a morphism in $^{B}_{H} \Bbb L$.
In fact, for any $m\in M, n\in N$ and $h\in H$, we have
\begin{eqnarray*}
&&C_{M,N}(h_{1}\c m\o h_{2}\c n)\\
&=&\langle(h_{1}\c m)_{(-1)}|(h_{2}\c n)_{(-1)}\rangle
         R^{(2)}\c \nu^{-2}(h_{2}\c n)_{(0)}\o R^{(1)}\c \mu^{-2}(h_{1}\c m)_{(0)}\\
&\stackrel{(2.1)}{=}&\langle\beta(m_{(-1)})|\beta(n_{(-1)})\rangle
         R^{(2)}\c \nu^{-2}(\alpha(h_{2})\c n_{(0)})\o R^{(1)}\c \mu^{-2}(\alpha(h_{1})\c m_{(0)})\\
&\stackrel{(HM2)}{=}&\langle m_{(-1)}|n_{(-1)}\rangle
        \alpha^{-1}(R^{(2)}h_{2})\c \nu^{-1}(n_{(0)})\o\alpha^{-1}(R^{(1)}h_{1})\c \mu^{-1}(m_{(0)}),\\
&&h\c C_{M,N}(m\o n)\\
&=&\langle m_{(-1)}|n_{(-1)}\rangle
         h\c(R^{(2)}\c \nu^{-2}(n_{(0)})\o R^{(1)}\c\mu^{-2}( m_{(0)}))\\
&=&\langle m_{(-1)}|n_{(-1)}\rangle
        h_{1}\c(\a^{-1}(R^{(2)})\c \nu^{-2}(n_{(0)}))\o h_{2}\c(\a^{-1}(R^{(1)})\c \mu^{-2}(m_{(0)}))\\
&\stackrel{(HM2)}{=}&\langle m_{(-1)}|n_{(-1)}\rangle
       \alpha^{-1}(h_{1}R^{(2)})\c \nu^{-1}(n_{(0)})\o \alpha^{-1}(h_{2}R^{(1)})\c \mu^{-1}(m_{(0)})\\
&\stackrel{(QHA4)}{=}&\langle m_{(-1)}|n_{(-1)}\rangle
        \alpha^{-1}(R^{(2)}h_{2})\c \nu^{-1}(n_{(0)})\o\alpha^{-1}(R^{(1)}h_{1})\c \mu^{-1}(m_{(0)}).
\end{eqnarray*}
The third equality holds since  $\langle |\rangle$ is $\beta$-invariant and the fifth equality holds since $R$ is $\alpha$-invariant.
So $C_{M,N}$ is left $(H,\alpha)$-linear.
Similarly, one may check that $C_{M,N}$ is left $(B,\beta)$-colinear.

Now we prove that the braiding  $C_{M,N}$ is natural.
For any $(M,\mu), (M',\mu'),$ $(N,\nu), (N',\nu')$ $\in{}^{B}_{H} \Bbb L$,
let $f: M\rightarrow M'$ and $g: N\rightarrow N'$ be two morpshisms in $^{B}_{H} \Bbb L$,
it is sufficient to verify the identity $(g\o f)\circ C_{M,N}=C_{M',N'}\circ(f\o g)$.
For this purpose, we take $m\in M, n\in N$ and do the following calculation:
\begin{eqnarray*}
(g\o f)\circ C_{M,N}(m\o n)
&=&\langle m_{(-1)}|n_{(-1)}\rangle(g\o f)(R^{(2)}\c \nu^{-2}(n_{(0)})\o R^{(1)}\c \mu^{-2}(m_{(0)}))\\
&=&\langle m_{(-1)}|n_{(-1)}\rangle g(R^{(2)}\c \nu^{-2}(n_{(0)}))\o f(R^{(1)}\c \mu^{-2}(m_{(0)}))\\
&=&\langle m_{(-1)}|n_{(-1)}\rangle R^{(2)}\c g(\nu^{-2}(n_{(0)}))\o R^{(1)}\c f(\mu^{-2}(m_{(0)})),\\
C_{M',N'}\circ(f\o g)(m\o n)
&=&C_{M',N'}(f(m)\o g(n))\\
&=&\langle f(m)_{(-1)}|g(n)_{(-1)}\rangle R^{(2)}\c \nu^{-2}(g(n)_{(0)})\o (R^{(1)}\c \mu^{-2}(f(m)_{(0)})\\
&=&\langle m_{(-1)}|n_{(-1)}\rangle R^{(2)}\c \nu^{-2}(g(n_{(0)}))\o R^{(1)}\c \mu^{-2}(f(m_{(0)}))\\
&=&\langle m_{(-1)}|n_{(-1)}\rangle R^{(2)}\c g(\nu^{-2}(n_{(0)}))\o R^{(1)}\c f(\mu^{-2}(m_{(0)})).
\end{eqnarray*}
The sixth equality holds since $f,g$ are left $(B,\beta)$-colinear.
So the braiding $C_{M,N}$ is natural, as needed.

Next, we will show that the braiding $C_{M,N}$ is an isomorphsim with inverse map
\begin{eqnarray*}
C^{-1}_{M,N}: N\o M\rightarrow M\o N ,
n\o m\rightarrow \langle S^{-1}(m_{(-1)})|n_{(-1)}\rangle S(R^{(1)})\c \mu^{-2}(m_{(0)})\o R^{(2)}\c \nu^{-2}(n_{(0)}).
\end{eqnarray*}
For any $m\in M, n\in N$, we have
\begin{eqnarray*}
&&C^{-1}_{M,N}\circ C_{M,N}(m\o n)\\
&=&\langle m_{(-1)}|n_{(-1)}\rangle C^{-1}_{M,N}(R^{(2)}\c \nu^{-2}(n_{(0)})\o R^{(1)}\c \mu^{-2}(m_{(0)}))\\
&=&\langle m_{(-1)}|n_{(-1)}\rangle\langle S^{-1}(\beta^{-1}(m_{(0)(-1)}))| \beta^{-1}(n_{(0)(-1)})\rangle\\
  &&~~~~~~S(r^{(1)})\c \mu^{-2}(\alpha(R^{(2)})\c \mu^{-2}(m_{(0)(0)}))\o r^{(2)}\c \nu^{-2}(\alpha(R^{(1)})\c \nu^{-2}(n_{(0)(0)}))\\
&\stackrel{(HCM2)}{=}&\langle \beta^{-1}(m_{(-1)1})|\beta^{-1}(n_{(-1)1})\rangle\langle S^{-1}(\beta^{-1}(m_{(-1)2}))| \beta^{-1}(n_{(-1)2})\rangle\\
   &&~~~~~~S(r^{(1)})\c (\alpha^{-1}(R^{(2)})\c \mu^{-3}(m_{(0)}))\o r^{(2)}\c (\alpha^{-1}(R^{(1)})\c \nu^{-3}(n_{(0)}))\\
&\stackrel{(HM2)}{=}&\langle m_{(-1)1}|n_{(-1)1}\rangle\langle S^{-1}(m_{(-1)2})| n_{(-1)2}\rangle\\
   &&~~~~~~\alpha^{-1}(S(r^{(1)})R^{(2)})\c \mu^{-2}(m_{(0)})\o \alpha^{-1}(r^{(2)}R^{(1)})\c \nu^{-2}(n_{(0)})\\
&\stackrel{(CHA1)}{=}&\langle S^{-1}(\beta^{-1}(m_{(-1)2}))\beta^{-1}(m_{(-1)1})|\b(n_{(-1)})\rangle
    1_{H}\c \mu^{-2}(m_{(0)})\o1_{H}\c \nu^{-2}(n_{(0)})\\
&=&\langle\beta^{-2}(S^{-1}(m_{(-1)2})m_{(-1)1})|n_{(-1)}\rangle
  1_{H}\c \mu^{-2}(m_{(0)})\o 1_{H}\c \nu^{-2}(n_{(0)})\\
&=&\langle \epsilon(m_{(-1)})1_{H}|n_{(-1)}\rangle
   \mu^{-1}(m_{(0)})\o  \nu^{-1}(n_{(0)})\\
&=&\epsilon(m_{(-1)})\epsilon(n_{(-1)}) \mu^{-1}(m_{(0)})\o  \nu^{-1}(n_{(0)})\\
&=&m\o n.
\end{eqnarray*}
The second equality holds since $\rho(R^{(2)}\c \nu^{-2}(n_{(0)}))=\beta^{-1}(n_{(0)(-1)})\o \alpha(R^{(2)})\c n_{(0)(0)}$
 and the fifth equality holds since $R^{-1}=S(r^{(1)})\o r^{(2)}$.

Now let us verify the hexagon axioms ($H_{1}, H_{2}$) from Section XIII. 1.1 of \cite{Kassel}.
We need to show that the following diagram ($H_{1}$) commutes
for any $(U,\mu), (V,\nu), (W,\omega)\in{}^{B}_{H} \Bbb L$:
$$\aligned \xymatrix{
(U\otimes V)\otimes W \ar[d]_{C_{U,V} \o id_W} \ar[rr]^{a_{U,V,W}} && U\otimes (V\otimes W) \ar[rr]^{C_{U,V \o W}} & & (V \otimes W)\otimes U \ar[d]^{a_{V,W,U}}  \\
(V \o U) \o W \ar[rr]^{a_{V,U,W}} & & V \o (U \o W) \ar[rr]^{id_V \o C_{U,W}} && V \o (W \o U),}
 \endaligned$$

For this purpose, let $u\in U, v\in V, w\in W$, then we have
\begin{eqnarray*}
&&a_{V,U,W}\circ C_{U,V\o W}\circ a_{U,V,W}((u\o v)\o w)\\
&=&a_{V,U,W}\circ C_{U,V\o W}(\mu^{-1}(u)\o(v\o \omega(w)))\\
&=&\langle \beta^{-1}(u_{(-1)})|\b^{-2}(v_{(-1)})\beta^{-1}(w_{(-1)})\rangle a_{V,U,W}\\
   &&~~~~~~~~~~(R^{(2)}\c (\nu^{-2}\o\omega^{-2})(v_{(0)}\o \omega(w_{(0)}))\o R^{(1)}\c \mu^{-3}(u_{(0)}))\\
&=&\langle \beta(u_{(-1)})|v_{(-1)}\beta(w_{(-1)})\rangle a_{V,U,W}\\
   &&~~~~~~~~~~(R^{(2)}\c (\nu^{-2}(v_{(0)})\o \omega^{-1}(w_{(0)}))\o R^{(1)}\c \mu^{-3}(u_{(0)}))\\
&=&\langle \beta(u_{(-1)})|v_{(-1)}\beta(w_{(-1)})\rangle\\
   &&~~~~~~~~~~\a^{-1}(R^{(2)}_{1})\c \nu^{-3}(v_{(0)})\o(R^{(2)}_{2}\c \omega^{-1}(w_{(0)})\o\a(R^{(1)})\c \mu^{-2}(u_{(0)}))\\
&\stackrel{(QHA3)}{=}&\langle \beta(u_{(-1)})|v_{(-1)}\beta(w_{(-1)})\rangle\\
    &&~~~~~~~~~~r^{(2)}\c \nu^{-3}(v_{(0)})\o(\a(R^{(2)})\c \omega^{-1}(w_{(0)})\o(R^{(1)}r^{(1)})\c  \mu^{-2}(u_{(0)}))
\end{eqnarray*}
and
\begin{eqnarray*}
&&(id_{V}\o C_{U,W})\circ a_{V,U,W}\circ(C_{U,V}\o id_{W})((u\o v)\o w)\\
&=&\langle u_{(-1)}|v_{(-1)}\rangle(id_{V}\o C_{U,W})\circ a_{V,U,W}((R^{(2)}\c \nu^{-2}(v_{(0)})\o R^{(1)}\c \mu^{-2}(u_{(0)}))\o w)\\
&=&\langle u_{(-1)}|v_{(-1)}\rangle(id_{V}\o C_{U,W})
   \a^{-1}(R^{(2)})\c \nu^{-3}(v_{(0)})\o(R^{(1)}\c \mu^{-2}(u_{(0)})\o \omega(w))\\
&=&\langle u_{(-1)}|v_{(-1)}\rangle\langle \beta^{-1}(u_{(0)(-1)})|\beta(w_{(-1)})\rangle\\
  &&~~~~~~\a^{-1}(R^{(2)})\c \nu^{-3}(v_{(0)})\o(r^{(2)}\c\omega^{-1}(w_{(0)})\o r^{(1)}\c\mu^{-2}(\a(R^{(1)})\c \mu^{-2}(u_{(0)(0)})))\\
&\stackrel{(HCM2)}{=}&\langle \beta^{-1}(u_{(-1)1})|v_{(-1)}\rangle\langle \beta^{-1}(u_{(-1)2})|\beta(w_{(-1)})\rangle\\
   &&~~~~~~\a^{-1}(R^{(2)})\c \nu^{-3}(v_{(0)})\o(r^{(2)}\c\omega^{-1}(w_{(0)})\o\a^{-1}(r^{(1)}R^{(1)})\c \mu^{-2}(u_{(0)}))\\
&\stackrel{(CHA2)}{=}&\langle u_{(-1)}|\b^{-1}(v_{(-1)})w_{(-1)}\rangle\\
  &&~~~~~~\a^{-1}(R^{(2)})\c \nu^{-3}(v_{(0)})\o(r^{(2)}\c\omega^{-1}(w_{(0)})\o\a^{-1}(r^{(1)}R^{(1)})\c \mu^{-2}(u_{(0)}))\\
&=&\langle \beta(u_{(-1)})|v_{(-1)}\beta(w_{(-1)})\rangle\\
  &&~~~~~~ R^{(2)}\c \nu^{-3}(v_{(0)})\o(\a(r^{(2)})\c \omega^{-1}(w_{(0)})\o(r^{(1)}R^{(1)})\c  \mu^{-2}(u_{(0)}))
\end{eqnarray*}
Since $r=R$, it follows that
$a_{V,U,W}\circ C_{U,V\o W}\circ a_{U,V,W}=(id_{V}\o C_{U,W})\circ a_{V,U,W}\circ(C_{U,V}\o id_{W})$,
that is, the diagram ($H_{1}$) commutes.

Now we check that the diagram ($H_{2}$) commutes for any $(U,\mu), (V,\nu), (W,\omega)\in{} ^{B}_{H}\Bbb{L}$:

$$\aligned
\xymatrix{
U\otimes (V\otimes W) \ar[d]_{id_U \o C_{V,W}} \ar[rr]^{a^{-1}_{U,V,W}} && (U\otimes V)\otimes W \ar[rr]^{C_{U\o V, W}} & & W \otimes (U\otimes V) \ar[d]^{a^{-1}_{W,U,V}}  \\
U \o (W\o V) \ar[rr]^{a^{-1}_{U,W, V}} & & (U \o W)\o V \ar[rr]^{C_{U,W} \o id_V } && (W \o U) \o V .}
\endaligned$$

In fact, for any  $u\in U, v\in V, w\in W$, we obtain
\begin{eqnarray*}
&&a^{-1}_{W,U,V}\circ C_{U\o V, W}\circ a^{-1}_{U,V,W}(u\o (v\o w))\\
&=&a^{-1}_{W,U,V}\circ C_{U\o V, W}((\mu(u)\o v)\o \omega^{-1}(w))\\
&=&\langle \beta^{-1}(u_{(-1)})\beta^{-1}(v_{(-2)})|\beta^{-1}(w_{(-1)})\rangle a^{-1}_{W,U,V}\\
     &&~~~~~~~~~~(R^{(2)}\c \omega^{-3}(w_{(0)})\o R^{(1)}\c (\mu^{-1}(u_{(0)})\o \nu^{-2}(v_{(0)})))\\
&=&\langle \beta(u_{(-1)})v_{(-1)}|\beta(w_{(-1)})\rangle a^{-1}_{W,U,V}\\
     &&~~~~~~~~~~(R^{(2)}\c \omega^{-3}(w_{(0)})\o (R^{(1)}_{1}\c \mu^{-1}(u_{(0)})\o R^{(1)}_{2}\c  \nu^{-2}(v_{(0)})))\\
&=&\langle \beta(u_{(-1)})v_{(-1)}|\beta(w_{(-1)})\rangle\\
    &&~~~~~~~~~~(\omega(R^{(2)}\c \omega^{-2}(w_{(0)}))\o R^{(1)}_{1}\c \mu^{-1}(u_{(0)}))\o\a^{-1}(R^{(1)}_{2})\c \nu^{-3}(v_{(0)})\\
&=&\langle \beta(u_{(-1)})v_{(-1)}|\beta(w_{(-1)})\rangle\\
   &&~~~~~~~~~~ (\alpha^{-1}(R^{(2)})\c\omega^{-2}(w_{(0)})\o R^{(1)}_{1}\c \mu^{-1}(u_{(0)}))\o \alpha(R^{(1)}_{2})\c \nu(v_{(0)})\\
&\stackrel{(QHA2)}{=}&\langle \beta(u_{(-1)})v_{(-1)}|\beta(w_{(-1)})\rangle\\
    &&~~~~~~~~~~(\alpha^{-1}(R^{(2)}r^{(2)})\c\omega^{-2}(w_{(0)})\o R^{(1)}\c \mu^{-1}(u_{(0)}))\o \alpha^{-1}(r^{(1)})\c \nu^{-3}(v_{(0)}).
\end{eqnarray*}
Also we can get
\begin{eqnarray*}
&&(C_{U,W}\o id_{V})\circ a^{-1}_{U,W,V}\circ(id_{U}\o C_{V,W})(u\o (v\o w))\\
&=&\langle v_{(-1)})|w_{(-1)}\rangle(C_{U,W}\o id_{V})\circ a^{-1}_{U,W,V}
    (u\o (R^{(2)}\c \omega^{-2}(w_{(0)})\o R^{(1)}\c \nu^{-2}(v_{(0)})))\\
&=&\langle v_{(-1)})|w_{(-1)}\rangle(C_{U,W}\o id_{V})
    ((\mu(u)\o R^{(2)}\c\omega^{-2}(w_{(0)}))\o\a^{-1}(R^{(1)})\c \nu^{-3}(v_{(0)}))\\
&=&\langle v_{(-1)})|w_{(-1)}\rangle\langle \beta(u_{(-1)})|\beta^{-1}(w_{(0)(-1)})\rangle\\
   &&~~~~~(r^{(2)}\c\omega^{-2}(\alpha(R^{(2)})\c \omega^{-2}(w_{(0)(0)}))\o  r^{(1)}\c \mu^{-1}(u_{(0)}))\o\a^{-1}(R^{(1)})\c \nu^{-3}(v_{(0)})\\
&\stackrel{(HCM2)}{=}&\langle v_{(-1)})|\beta^{-1}(w_{(-1)1})\rangle\langle \beta(u_{(-1)})|\beta^{-1}(w_{(-1)2})\rangle\\
   &&~~~~~(r^{(2)}\c(\a^{-1}(R^{(2)})\c \omega^{-3}(w_{(0)}))\o  r^{(1)}\c \mu^{-1}(u_{(0)}))\o\alpha^{-1}(R^{(1)})\c \nu^{-3}(v_{(0)})\\
&\stackrel{(CHA1)}{=}&\langle u_{(-1)}\beta^{-1}(v_{(-1)})|w_{(-1)}\rangle\\
     &&~~~~~(\alpha^{-1}(r^{(2)}R^{(2)})\c \omega^{-2}(w_{(0)})\o  r^{(1)}\c \mu^{-1}(u_{(0)}))\o\alpha^{-1}(R^{(1)})\c \nu^{-3}(v_{(0)}).
\end{eqnarray*}
So the diagram ($H_{2}$) commutes since  $r=R$.
This ends the proof.
\medskip

\noindent{\bf Corollary  3.2.}
Under the hypotheses of the Theorem 3.1, the braiding $C$ is a solution of the quantum Yang-Baxter equation
\begin{eqnarray*}
&&(id_{W}\o C_{U,V})\circ a_{W,U,V}\circ( C_{U,W}\o id_{V})\circ a^{-1}_{W,V,U}\circ( id_{U}\o C_{V,W})\circ a_{U,V,W}\nonumber\\
&=&a_{W,V,U}\circ( C_{W,V}\o id_{U})\circ a^{-1}_{W,V,U}\circ( id_{V}\o C_{U,W})\circ a_{V,U,W}\circ( C_{U,V}\o id_{W}).
\end{eqnarray*}

\noindent{\bf Proof.} Straightforward.
\medskip

\noindent{\bf Lemma 3.3.}
Let $(H,R,\alpha)$ be a quasitriangular Hom-Hopf algebra
and $(B,\langle|\rangle,\beta)$ a coquasitriangular Hom-Hopf algebra.
Define a linear map
\begin{eqnarray*}
(H\o B)\o M\rightarrow M,(h\o x)\rightharpoonup m=\langle x|m_{(-1)}\rangle \a^{-3}(h)\c \mu^{-1}(m_{(0)}),
\end{eqnarray*}
for any $h\in H, x\in B$ and $m\in(M,\mu)\in{}^{B}_{H} \Bbb L$.
Then $(M,\mu)$ becomes a left $(H\o B)$-Hom-module.
 \medskip

\noindent{\bf Proof.}
It is sufficient to show that the Hom-module action defined above satisfies  Definition 1.2.
For any $h,g\in H, x,y\in B$ and $m\in M$, we have
\begin{eqnarray*}
(1_{H}\o 1_{B})\rightharpoonup m
=\langle 1_{B}|m_{(-1)}\rangle 1_{H}\c \mu^{-1}(m_{(0)})
=\epsilon(m_{(-1)})m_{(0)}=\mu(m).
\end{eqnarray*}
That is, $(1_{H}\o 1_{B})\rightharpoonup m=\mu(m)$.
For the equality $\mu((h\o x)\rightharpoonup m)=(\alpha(h)\o \beta(x))\rightharpoonup \mu(m)$, we have
\begin{eqnarray*}
(\alpha(h)\o \beta(x))\rightharpoonup \mu(m)
&=&\langle \beta(x)|\beta(m_{(-1)})\rangle \alpha^{-2}(h)\c m_{(0)}\\
&=&\langle x|m_{(-1)}\rangle \alpha^{-2}(h)\c m_{(0)}
=\mu((h\o x)\rightharpoonup m),
\end{eqnarray*}
as required. Finally, we check the expression
$((h\o x)(g\o y))\rightharpoonup \mu(m)=(\alpha(h)\o\beta(x))\rightharpoonup((g\o y)\rightharpoonup m)$.
For this, we calculate
\begin{eqnarray*}
&&(\alpha(h)\o\beta(x))\rightharpoonup ((g\o y)\rightharpoonup m)\\
&=&\langle y|m_{(-1)}\rangle(\alpha(h)\o\beta(x))\c (\alpha^{-3}(g)\c \mu^{-1}(m_{(0)}))\\
&=&\langle y|m_{(-1)}\rangle\langle \beta(x)|m_{(0)(-1)}\rangle\alpha^{-2}(h)\c(\alpha^{-3}(g)\c\mu^{-2}(m_{(0)(0)}))\\
&\stackrel{(HCM2)}{=}&\langle y|\beta^{-1}(m_{(-1)1})\rangle\langle x|\beta^{-1}(m_{(-1)2})\rangle \alpha^{-3}(hg)\c m_{(0)}\\
&\stackrel{(CHA1)}{=}&\langle xy|\beta(m_{(-1)})\rangle\alpha^{-3}(hg)\c m_{(0)}\\
&=&((h\o x)(g\o y))\rightharpoonup \mu(m).
\end{eqnarray*}
So $(M,\mu)$ is a left $(H\o B)$-Hom-module. The proof is completed.
 \medskip

 \noindent{\bf Lemma 3.4.}
Let $(H,R,\alpha)$ be a quasitriangular Hom-Hopf algebra
and $(B,\langle|\rangle,\beta)$ a coquasitriangular Hom-Hopf algebra.
Define a linear map
\begin{eqnarray*}
\overline{\rho}: M\rightarrow (H\o B)\o M,~\overline{\rho}(m)=m_{[-1]}\o m_{[0]}= R^{(2)} \o \beta^{-3}(m_{(-1)})\o R^{(1)}\c \mu^{-1}(m_{(0)}),
\end{eqnarray*}
for any $m\in (M,\mu)$.
Then $(M,\mu)$ becomes a left $(H\o B)$-Hom-comodule.
 \medskip

\noindent{\bf Proof.}
We first show that $\overline{\rho}$ satisfies Eq. (HCM2). On the one side, we have
\begin{eqnarray*}
&&\Delta(m_{[-1]})\o \mu(m_{[0]})\\
&=&( R^{(2)}_{1} \o\beta^{-3}(m_{(-1)1}))\o( R^{(2)}_{2} \o\beta^{-3}(m_{(-1)2}))\o\alpha(R^{(1)})\c m_{(0)}\\
&=&(\a(r^{(2)})\o\beta^{-2}(m_{(-1)}))\o(\a(R^{(2)})\o\beta^{-3}(m_{(0)(-1)}))\o\alpha(R^{(1)})(r^{(1)}\c \mu^{-2}(m_{(0)(0)})).
\end{eqnarray*}
On the other side, we have
\begin{eqnarray*}
&&~~(\alpha\o\beta)(m_{[-1]})\o \overline{\rho}(m_{[0]})\\
&&=(\a(r^{(2)})\o \b^{-2}(m_{(-1)}))\o(R^{(2)} \o \beta^{-3}( ( r^{(1)} \c \mu^{-1}(m_{(0)}) )_{(-1)} ) \o R^{(1)} \\
&&~~~~~~~~~~\c \mu^{-1}(( r^{(1)} \c \mu^{-1}(m_{(0)}) )_{(0)} )   \\
&&= (\a(r^{(2)})\o\beta^{-2}(m_{(-1)}))\o(R^{(2)}\o\beta^{-3}(m_{(0)(-1)}))\o R^{(1)} \c (r^{(1)}\c \mu^{-2}(m_{(0)(0)})).
\end{eqnarray*}
Since $R$ is $\alpha$-invariant, we have $\Delta(m_{[-1]})\o \mu(m_{[0]})=(\alpha\o\beta)(m_{[-1]})\o \overline{\rho}(m_{[0]})$, as needed.

For Eq. (HCM1), we have
\begin{eqnarray*}
(\epsilon_{H}\o\epsilon_{B})(m_{[-1]})m_{[0]}
&=&\epsilon_{H}(R^{(2)})\epsilon_{B}(m_{(-1)})R^{(1)}\c \mu^{-1}(m_{(0)})\\
&=&1_{H}\c m
=\mu(m),\\
(\alpha\o\beta)(m_{[-1]})\o \mu(m_{[0]})
&=&(\alpha(R^{(2)})\o \beta^{-2}(m_{(-1)}))\o \mu(R^{(1)}\c \mu^{-1}(m_{(0)}))\\
&=& R^{(2)} \o \beta^{-3}(\beta(m_{(-1)}))\o R^{(1)}\c \mu^{-1}(\mu(m_{(0)}))\\
&=& \overline{\rho}(\mu(m)),
\end{eqnarray*}
as desired. And this finishes the proof.
\medskip

\noindent{\bf Theorem 3.5.}
Let $(H,R,\alpha)$ be a quasitriangular Hom-Hopf algebra
and $(B,\langle|\rangle,\beta)$ a coquasitriangular Hom-Hopf algebra.
Then the Hom-Long dimodules category $^{B}_{H} \Bbb L$ is a monoidal subcategory
 of Hom-Yetter-Drinfeld category $^{H\o B}_{H\o B}\Bbb {YD}$.
 \medskip

\noindent{\bf Proof.}
Let $m\in(M,\mu)\in {}^{B}_{H}\mathcal{L}$ and $h\in H$.
Here we first note that
$\rho(h\c \mu^{-1}(m_{(0)}))=m_{(0)(-1)}\o\alpha(h)\c \mu^{-1}(m_{(0)(0)})$.
It is sufficient to show that the left $(H\o B)$-Hom-module action in Lemma 3.3 and
the left $(H\o B)$-Hom-comodule structure in Lemma 3.4 satisfy the compatible condition Eq. (HYD).
Indeed, for any $h \in H$, $x \in B$, $m \in M$, we have
\begin{eqnarray*}
&&(h_1 \o x_1)(\a \o \b )(m_{[-1]}) \o (\a^3(h_2) \o \b^3(x_2))\rightharpoonup m_{[0]} \\
&=& h_1 \a(R^{(2)}) \o x_1 \b^{-2}(m_{(-1)}) \o \langle \b^3(x_2) | (R^{(1)} \c \mu^{-1}(m_{(0)}))_{(-1)} \rangle h_2 \c \mu^{-1}( (R^{(1)} \c \mu^{-1}(m_{(0)}))_{(0)} )\\
&=& h_1 \a(R^{(2)}) \o x_1 \b^{-3}(m_{(-1)1}) \o \langle \b^3(x_2) | m_{(-1)2} \rangle h_2 \c ( R^{(1)} \c \mu^{-1}(m_{(0)}) )\\
&=& h_1 \a(R^{(2)}) \o x_1 \b^{-3}(m_{(-1)1}) \o \langle x_2 | \b^{-3}(m_{(-1)2}) \rangle \a^{-1}(h_2 \a(R^{(1)}) ) \c m_{(0)} \\
&=& R^{(2)} h_2 \o \b^{-3}(m_{(-1)2}) x_2 \langle x_1 | \b^{-3}(m_{(-1)1}) \rangle  \o (\a^{-1}(R^{(1)}) \a^{-1}(h_1)) \c m_{(0)} \\
&=& \langle \a^2(x_1) | m_{(-1)} \rangle R^{(2)} h_2 \o \b^{-3}(m_{(0)(-1)}) x_2 \o (\a^{-1}(R^{(1)}) \a^{-1}(h_1)) \c \mu^{-1}(m_{(0)(0)}) \\
&=& \langle \a^2(x_1) | m_{(-1)} \rangle (R^{(2)}  \o \b^{-3}( {\a^{-1}(h_1) \c \mu^{-1}(m_{(0)}) }_{(-1)} ) )(h_2 \o x_2) \\
&&~~~~~~~~~~~~\o R^{(1)} \c \mu^{-1}({\a^{-1}(h_1) \c \mu^{-1}(m_{(0)}) }_{(0)})   \\
&=& {(\a^2(h_1) \o \b^2(x_1))\rightharpoonup m}_{[-1]} (h_2 \o x_2)  \o {(\a^2(h_1) \o \b^2(x_1))\rightharpoonup m}_{[0]}.
\end{eqnarray*}
So $(M,\mu)\in{}^{H\o B}_{H\o B}\Bbb{HYD}$. The proof is completed.
 \medskip

\noindent{\bf Proposition 3.6.}
Under the hypotheses of the Theorem 3.5,
 $^{B}_{H}\Bbb{L}$ is a braided monoidal subcategory
 of  $^{H\o B}_{H\o B}\Bbb{HYD}$.
 \medskip

\noindent{\bf Proof.}
It is sufficient to show that the braiding in the category
$^{B}_{H}\Bbb{L}$ is compatible to the braiding in $^{H\o B}_{H\o B}\Bbb{HYD}$.
In fact, for any $m\in (M,\mu)$ and $n\in (N,\nu)$, we have
\begin{eqnarray*}
C_{M,N}(m\o n)
&=& (\a^2(R^{(2)}) \o \b^{-1}(m_{(-1)}))\rightharpoonup \nu^{(-1)}(n) \o \a^{-1}(R^{(1)}) \c \mu^{-2}(m_{(0)}) \\
&=& \langle \b^{-1}(m_{(-1)}) | \b^{-1}(n_{(-1)}) \rangle \a^{-1}(R^{(2)}) \c \nu^{-2}(n_{(0)}) \o \a^{-1}(R^{(1)}) \c \mu^{-2}(m_{(0)}) \\
&=& \langle m_{(-1)} |  n_{(-1)} \rangle R^{(2)} \c \nu^{-2}(n_{(0)}) \o R^{(1)} \c \mu^{-2}(m_{(0)}),
\end{eqnarray*}
as desired.This finishes the proof.

\section{Symmetries in Hom-Long dimodule categories}
\def\theequation{\arabic{section}.\arabic{equation}}
\setcounter{equation} {0}

In this section, we obtain a sufficient condition for the Hom-Long dimodule category $^{B}_{H} \Bbb L$ to be symmetric.
\medskip

Let $\mathcal{C}$ be a monoidal category and $C$ a braiding on $\mathcal{C}$.
The braiding $C$ is called a symmetry \cite{Joyal1993, Kassel} if
$C_{Y,X}\circ C_{X,Y}=id_{X\otimes Y}$ for all $X,Y\in \mathcal{C}$,
and the category $\mathcal{C}$ is called symmetric.
\medskip

\noindent{\bf Proposition 4.1.}
Let $(H,R,\alpha)$  be a triangular Hom-Hopf algebra  and
 $(B,\beta)$ a   Hom-Hopf algebra.
 Then the category  $_{H} \Bbb M$   of left  $(H,\alpha)$-Hom-modules is a
symmetric subcategory of $^{B}_{H} \Bbb L$  under the left $(B,\beta)$-comodule structure
$\rho(m)=1_{B}\otimes \mu(m)$,
 where $m\in (M,\mu)\in{}_{H} \Bbb M$,
 and the braiding is defined as
\begin{eqnarray*}
C_{M,N}: M\o N\rightarrow N\o M, m\o n\rightarrow R^{(2)}\c \nu^{-1}(n)\o R^{(1)}\c \mu^{-1}(m),
\end{eqnarray*}
for all $m\in (M,\mu)\in{}_{H} \Bbb M, n\in (N,\nu)\in{}_{H} \Bbb M.$
\medskip

\noindent{\bf Proof.}
It is clear that  $(M,\rho,\mu)$ is a left $(B,\beta)$-Hom-comodule under the left $(B,\beta)$-comodule structure given above.
Now we check that the left $(B,\beta)$-comodule structure satisfies the compatible condition Eq. (2.1).
For this purpose, we take $h\in H, m\in(M,\mu)\in{}_{H} \Bbb M$, and calculate
\begin{eqnarray*}
\rho(h\c m)
=1_{B}\otimes \mu(h\c m)
=1_{B}\otimes \alpha(h)\c\mu(m)
=\beta(m_{(-1)})\otimes \alpha(h)\c m_{(0)}.
\end{eqnarray*}
So, Eq. (2.1) holds.  That is, $(M,\rho,\mu)$ is an $(H,B)$-Hom-Long dimodule.

Next we verify that any morphism in $_{H} \Bbb M$ is left $(B,\beta)$-colinear, too.
Indeed,  for any $m\in (M,\mu)\in{}_{H} \Bbb M$ and $n\in (N,\nu)\in{}_{H} \Bbb M$.
Assume that $f: (M,\mu)\rightarrow (N,\nu)$ is a morphism in $_{H} \Bbb M$,  then
\begin{eqnarray*}
(id_{B}\otimes f)\rho(m)
=1_{B}\otimes f(\mu(m))
=1_{B}\otimes\nu(f(m))
=\rho(f(m)).
\end{eqnarray*}
So $f$ is left $(B,\beta)$-colinear, as desired.
Therefore, $_{H} \Bbb M$  is a subcategory of $^{B}_{H} \Bbb L$.

Finally, we prove that $_{H} \Bbb M$  is a symmetric subcategory of $^{B}_{H} \Bbb L$.
Since
$
C_{M,N}(m\o n)
= R^{(2)}\c \nu^{-1}(n)\o R^{(1)}\c \mu^{-1}(m),
$
 for all $m\in (M,\mu)\in{}_{H} \Bbb M$ and $n\in (N,\nu)\in{}_{H} \Bbb M$,
 we have
\begin{eqnarray*}
C_{N,M}\circ C_{M,N}(m\o n)
&=&C_{N,M}(R^{(2)}\c \nu^{-1}(n)\o R^{(1)}\c \mu^{-1}(m))\\
&=&r^{(2)}\c\mu^{-1}(R^{(1)}\c \mu^{-1}(m))\o r^{(1)}\c\nu^{-1}(R^{(2)}\c \nu^{-1}(n))\\
&=&r^{(2)}\c(\a^{-1}(R^{(1)})\c \mu^{-2}(m))\o r^{(1)}\c(\a^{-1}(R^{(2)})\c \nu^{-2}(n))\\
&=&\alpha^{-1}(r^{(2)}R^{(1)})\c \mu^{-1}(m)\o \alpha^{-1}(r^{(1)}R^{(2)})\c \nu^{-1}(n)\\
&=&1_{H}\c \mu^{-1}(m)\o 1_{H}\c \nu^{-1}(n)
=m\o n.
\end{eqnarray*}
It follows that the braiding $C_{M,N}$ is symmetric.
The proof is completed.
\medskip

\noindent{\bf Proposition 4.2.}
Let $(B,\langle|\rangle,\beta)$ be a cotriangular Hom-Hopf algebra  and
 $(H,\alpha)$ a Hom-Hopf algebra.
 Then the category $^{B} \Bbb M$ of left  $(B,\beta)$-Hom-comodules is a
symmetric subcategory of $^{B}_{H} \Bbb L$ under the left $(H,\alpha)$-module action
$h\c m=\epsilon(h)\mu(m)$,
 where $h\in H, m\in (M,\mu)\in{}^{B} \Bbb M$,
 and the braiding is given by
\begin{eqnarray*}
C_{M,N}: M\o N\rightarrow N\o M, m\o n\rightarrow\langle m_{(-1)}|n_{(-1)}\rangle \nu^{-2}(n_{(0)})\o \mu^{-2}(m_{(0)}),
\end{eqnarray*}
for all $m\in (M,\mu)\in{}^{B} \Bbb M, n\in (N,\nu)\in{}^{B} \Bbb M.$
\medskip

\noindent{\bf Proof.}
We first show that the left $(H,\alpha)$-module action defined above forces
 $(M,\mu)$ to be a left $(H,\alpha)$-module, but this is easy to check.
For the compatible condition Eq. (2.1), we take $h\in H, m\in(M,\mu)\in{}^{B} \Bbb M$
 and calculate as follows:
\begin{eqnarray*}
\rho(h\c m)=1_{B}\o\mu(h\c m)=1_{B}\o\epsilon(h)\mu( m)=\beta(m_{(-1)})\otimes \alpha(h)\c m_{(0)}.
\end{eqnarray*}
So, Eq. (2.1) holds, as required.  Therefore, $(M,\rho,\mu)$ is an $(H,B)$-Hom-Long dimodule.

Now we verify that any morphism in $^{B} \Bbb M$ is left $(H,\alpha)$-linear, too.
Indeed,  for any $m\in (M,\mu)\in{}^{B} \Bbb M$ and $n\in (N,\nu)\in{}^{B} \Bbb M$.
Assume that $f: (M,\mu)\rightarrow (N,\nu)$ is a morphism in $^{B} \Bbb M$,  then
\begin{eqnarray*}
f(h\c m)
=f(\epsilon(h)\mu(m))
=\epsilon(h)\mu(f(m))
=h\c f(m).
\end{eqnarray*}
So $f$ is left $(H,\alpha)$-linear, as desired.
Therefore, $^{B} \Bbb M$  is a subcategory of $^{B}_{H} \Bbb L$.

Finally, we show that  $^{B} \Bbb M$ is a symmetric subcategory of $^{B}_{H} \Bbb L$.
Since
$
C_{M,N}(m\o n)
=\langle m_{(-1)}|n_{(-1)}\rangle \nu^{-1}(n_{(0)})\o \mu^{-1}(m_{(0)}),
$
 for all $m\in (M,\mu)\in{}^{B} \Bbb M$ and $n\in (N,\nu)\in{}^{B} \Bbb M$,
then
\begin{eqnarray*}
&&C_{N,M}\circ C_{M,N}(m\o n)\\
&=&\langle m_{(-1)}|n_{(-1)}\rangle C_{N,M}(\nu^{-1}(n_{(0)})\o \mu^{-1}(m_{(0)}))\\
&=&\langle m_{(-1)}|n_{(-1)}\rangle\langle \b^{-1}(n_{(0)(-1)})|\b^{-1}(m_{0(-1)})\rangle(\mu^{-2}(m_{(0)(0)})\o \nu^{-2}(n_{(0)(0)})\\
&=&\langle \beta^{-1}(m_{(-1)1})|\beta^{-1}(n_{(-1)1})\rangle\langle \b^{-1}(n_{(-1)2})|\b^{-1}(m_{(-1)2})\rangle\mu^{-1}(m_{(0)})\o \nu^{-1}(n_{(0)})\\
&=&\epsilon(m_{(-1)})\epsilon(n_{(-1)})\mu^{-1}(m_{(0)})\o \nu^{-1}(n_{(0)})
=m\o n,
\end{eqnarray*}
where the fourth equality holds since $\langle | \rangle$ is $\b$-invariant.
It follows that the braiding $C_{M,N}$ is symmetric.
The proof is completed.
\medskip

\noindent{\bf Theorem 4.3.}
Let $(H,\alpha)$  be a triangular Hom-Hopf algebra  and $(B,\langle|\rangle,\beta)$  a cotriangular Hom-Hopf algebra.
 Then the category $^{B}_{H} \Bbb L$ is symmetric.
\medskip

\noindent{\bf Proof.}
For any $m\in (M,\mu)\in{}^{B}_{H} \Bbb L$ and $n\in (N,\nu)\in{}^{B}_{H} \Bbb L$, we have
\begin{eqnarray*}
&&C_{N,M}\circ C_{M,N}(m\o n)\\
&=&\langle m_{(-1)}|n_{(-1)}\rangle C_{N,M}(R^{(2)}\c \nu^{-2}(n_{(0)})\o R^{(1)}\c \mu^{-2}(m_{(0)}))\\
&=&\langle m_{(-1)}|n_{(-1)}\rangle\langle \beta(n_{(0)(-1)})|\beta(m_{(0)(-1)})\rangle \\
&&~~~~~~~~r^{(2)}\c\mu^{-2}(\alpha(R^{(1)})\c\mu^{-2}( m_{(0)(0)}))\o r^{(1)}\c\nu^{-2}(\alpha(R^{(2)})\c \nu^{-2}(n_{(0)(0)}))\\
&=&\langle \beta^{-1}(m_{(-1)1})|\beta^{-1}(n_{(-1)1})\rangle\langle \beta^{-1}(n_{(-1)2})|\beta^{-1}(m_{(-1)2})\rangle\\
&&~~~~~~~~\alpha^{-1}(r^{(2)}R^{(1)})\c \mu^{-2}(m_{(0)})\o\alpha^{-1}(r^{(1)}R^{(2)})\c \nu^{-2}(n_{(0)})\\
&=&\epsilon(m_{(-1)})\epsilon(n_{(-1)}) 1_{H}\c \mu^{-2}(m_{(0)})\o 1_{H}\c \nu^{-2}(n_{(0)})\\
&=&\epsilon(m_{(-1)})\epsilon(n_{(-1)}) \mu^{-1}(m_{(0)})\o\nu^{-1}(n_{(0)})\\
&=&m\o n,
\end{eqnarray*}
as desired. This finishes the proof.

\section{New solutions of the Hom-Long Equation}
\def\theequation{\arabic{section}.\arabic{equation}}
\setcounter{equation} {0}

In this section, we will present a kind of new solutions of  the Hom-Long equation.
\medskip

\noindent{\bf Definition 5.1.}
Let $(H,\alpha)$  be a  Hom-bialgebra  and $(M,\mu)$  a   Hom-module over $(H,\a)$.
 Then $R\in End(M\o M)$ is called the solution of the Hom-Long equation if it satisfies the nonlinear equation:
 \begin{eqnarray}
R^{12}\circ R^{23}=R^{23}\circ R^{12},
\end{eqnarray}
where $R^{12}=R\o\mu,R^{23}=\mu\o R$.
\medskip

\noindent{\bf Example 5.2.}
 If $R\in End(M\o M)$ is invertible, then it is easy to see that $R$ is a solution of the Hom-Long equation if and only if $R^{-1}$ is too.
\medskip

\noindent{\bf Example 5.3.}
Let $(M,\mu)$  an $(H,\a)$-Hom-module with a basis $\{m_1,m_2,\cdots,m_n\}$.
Assume that $\mu$ is given by $\mu(m_i)=a_im_i$, where $a_i\in k,~i=1,2,\cdots,n$.
Define a map
\begin{eqnarray*}
R:~M\o M\rightarrow M\o M,~~R(m_i\o m_j)=b_{ij}m_i\o m_j,
\end{eqnarray*}
 where $b_{ij}\in k,~i,j=1,2,,\cdots,n.$
Then  $R$  is a solution of the Hom-Long equation (5.1).
Furthermore, if $a_i=1$, for all $i=1,2,\cdots,n$, then $R$ is a solution of the classical Long equation.
\medskip

\noindent{\bf Proposition 5.4.}
Let $(M,\mu)$  an $(H,\a)$-Hom-module with a basis $\{m_1,m_2,\cdots,m_n\}$.
Assume that $R,S\in End(M\o M,\mu\o\mu^{-1})$ given by the matrix formula
\begin{eqnarray*}
R(m_k\o m_l)= x_{kl}^{ij}m_i\o \mu^{-1}(m_j),~~S(m_k\o m_l)= y_{kl}^{ij}m_i\o \mu^{-1}(m_j),
\end{eqnarray*}
and $\mu(m_l)=z_{l}^{i}m_i$, where $x_{kl}^{ij},y_{kl}^{ij},z_{l}^{i}\in k$.
Then  $S^{12}\circ R^{23}=R^{23}\circ S^{12}$ if and only if
 \begin{eqnarray*}
z_{u}^{i}x_{vw}^{jk}y_{ij}^{pq}=z_{i}^{p}x_{jw}^{qk}y_{uv}^{ij},
\end{eqnarray*}
for all $k,p,q,u,v,w=1,2,\cdots,n$.
In particular, $R$ is a solution of Hom-Long equation if and only if
 \begin{eqnarray*}
z_{u}^{i}x_{vw}^{jk}x_{ij}^{pq}=z_{i}^{p}x_{jw}^{qk}x_{uv}^{ij}.
\end{eqnarray*}

\noindent{\bf Proof.}
According to the definition of $R,S,\mu$, we have
\begin{eqnarray*}
S^{12}\circ R^{23}(m_{u}\o m_{v}\o m_{w})
&=&S^{12}(z_{u}^{i}m_{i}\o x_{vw}^{jk}m_{j}\o \mu^{-1}(m_{k}))\\
&=&z_{u}^{i}x_{vw}^{jk}y_{ij}^{pq}(m_{p}\o \mu^{-1}(m_{q})\o m_{k}),\\
R^{23}\circ S^{12}(m_{u}\o m_{v}\o m_{w})
&=&R^{23}(y_{uv}^{ij}m_{i}\o\mu^{-1}(m_{j})\o m_{w})\\
&=&y_{uv}^{ij}z_{i}^{p}x_{jw}^{qk}(m_{p}\o \mu^{-1}(m_{q})\o m_{k}).
\end{eqnarray*}
It follows that  $S^{12}\circ R^{23}=R^{23}\circ S^{12}$ if and only if
$
z_{u}^{i}x_{vw}^{jk}y_{ij}^{pq}=z_{i}^{p}x_{jw}^{qk}y_{uv}^{ij}.
$
Furthermore, $R^{12}\circ R^{23}=R^{23}\circ R^{12}$ if and only if
$
z_{u}^{i}x_{vw}^{jk}x_{ij}^{pq}=z_{i}^{p}x_{jw}^{qk}x_{uv}^{ij}.
$
The proof is completed.
\medskip

In the following proposition, we use the notation: for any $F\in End(M\o M)$, we denote $F^{12}=F\o \mu, F^{23}=\mu\o  F,
F^{13}=(id\o\tau )\circ (F \o \mu)\circ (id\o \tau)$, and $\tau^{(123)}(x\o y\o z)=(z,x,y).$
\smallskip

\noindent{\bf Proposition 5.5.}
Let $(M,\mu)$  an $(H,\a)$-Hom-module and $R\in End(M\o M)$. The following statements are equivalent:

(1) $R$ is a solution of the Hom-Long equation.

(2) $U=\tau\circ R$ is a solution of the equation:
$$U^{13}\circ U^{23}=\tau^{(123)}\circ U^{13}\circ U^{12}.$$

(3) $T=R\circ\tau$ is a solution of the  equation:
$$T^{12}\circ T^{13}=T^{23}\circ T^{13}\circ\tau^{(123)}.$$

(4) $W=\tau\circ R \circ \tau$ is a solution of the equation:
$$\tau^{(123)}\circ W^{23}\circ W^{13}=W^{12}\circ W^{13}\circ\tau^{(123)}.$$

\noindent{\bf Proof.}
We just prove$(1)\Leftrightarrow (2)$, and similar for $(1)\Leftrightarrow (3) $ and $(1)\Leftrightarrow (4).$
Since $R=\tau\circ U$,  $R$ is a solution of the Hom-Long equation if and only if  $R^{12}\circ R^{23}=R^{23}\circ R^{12}$,
that is,
\begin{eqnarray}
\tau^{12}\circ U^{12}\circ \tau^{23}\circ U^{23}=\tau^{23}\circ U^{23}\circ\tau^{12}\circ U^{12}.
\end{eqnarray}
While $\tau^{12}\circ U^{12}\circ \tau^{23}=\tau^{23}\circ\tau^{13}\circ U^{13}$ and
        $\tau^{23}\circ U^{23}\circ\tau^{12}=\tau^{23}\circ\tau^{12}\circ U^{13}$,
(5.2) is equivalent to
$$\tau^{23}\circ \tau^{13}\circ  U^{13}\circ U^{23}=\tau^{23}\circ\tau^{12}\circ U^{13}\circ U^{12},$$
which is equivalent to $U^{13}\circ U^{23}=\tau^{(123)}\circ U^{13}\circ U^{12}$ from the fact
$\tau^{23}\circ\tau^{12}=\tau^{(123)}$.
\medskip

Next we will present a new solution for Hom-Long equation by the Hom-Long dimodule structures.
For this, we give the notion of $(H,\a)$-Hom-Long dimodules.
\medskip

\noindent{\bf Definition 5.6.}
Let  $(H,\alpha)$ a  Hom-bialgebra.
A left-left $(H,\a)$-Hom-Long dimodule is a quadrupl  $(M,\c, \rho,\mu)$,
where $(M, \c, \mu)$ is a left  $(H,\alpha)$-Hom-module
and  $(M, \rho, \mu)$ is a left $(H,\alpha)$-Hom-comodule such that
\begin{eqnarray}
\rho(h\c m)=\a(m_{(-1)})\o \alpha(h)\c m_{0},
\end{eqnarray}
for all $h\in H$ and $m\in M$.
\medskip

\noindent{\bf Remark 5.7.}
Clearly, left-left $(H,\a)$-Hom-Long dimodules is a special case of $(H,B)$-Hom-Long dimodules in Definition 2.1 by setting $(H,\a)=(B,\b).$
\medskip

\noindent{\bf Example 5.8.}
Let  $(H,\a)$ be a Hom-bialgebra and  $(M,\c,\mu)$ be a left $(H,\a)$-Hom-module.
Define a left $(H,\a)$-Hom-module structure and a left $(H,\a)$-Hom-comodule structure on $(H\o M,\a\o \mu)$ as follows:
\begin{eqnarray*}
h\c(g\o m)=\a(g)\o h\c\mu(m), ~~\rho(g\o m)=g_1\o g_2\o\mu(m),
\end{eqnarray*}
for all $h,g\in H$ and $m\in M$.
Then $(H\o M,\a\o \mu)$ is an $(H,\a)$-Hom-Long dimodule.
\medskip

\noindent{\bf Example 5.9.}
Let  $(H,\a)$ be a Hom-bialgebra and   $(M,\rho,\mu)$ be a left $(H,\a)$-Hom-comodule.
Define a left $(H,\a)$-Hom-module structure and be a left $(H,\a)$-Hom-comodule structure on $(H\o M,\a\o \mu)$ as follows:
\begin{eqnarray*}
h\c(g\o m)=hg\o\mu(m), ~~\rho(g\o m)=m_{-1}\o \a(g)\o m_0,
\end{eqnarray*}
for all $h,g\in H$ and $m\in M$.
Then $(H\o M,\a\o \mu)$ is an $(H,\a)$-Hom-Long dimodule.
\medskip

\noindent{\bf Theorem 5.10.}
Let  $(H,\a)$ be a Hom-bialgebra and   $(M,\c,\rho,\mu)$ be a $(H,\a)$-Hom-Long dimodule.
Then the map
\begin{eqnarray}
R_{M}: M\o M\rightarrow M\o M,~~~~m\o n\mapsto n_{-1}\c m \o n_0,
\end{eqnarray}
is a solution of the Hom-Long equation, for any $m,n\in M.$
\medskip

\noindent{\bf Proof.} For any  $l,m,n\in M$, we calculate
\begin{eqnarray*}
R_{M}^{12}\circ R_{M}^{23}(l\o m\o n)
&=&R_{M}^{12}(\mu(l)\o n_{(-1)}\c m\o n_0)\\
&=&(n_{(-1)}\c m)_{(-1)}\c\mu(l)\o(n_{(-1)}\c m)_{0}\o \mu(n_0)\\
&=&\a(m_{(-1)})\c\mu(l)\o\a(n_{(-1)})\c m_{0}\o\mu(n_0),\\
R_{M}^{23}\circ R_{M}^{12}(l\o m\o n)
&=&R_{M}^{23}(m_{(-1)}\c l\o m_0\o \mu(n))\\
&=&\mu(m_{(-1)}\c l)\o\a(n_{(-1)}))\c m_0\o\mu(n_0)\\
&=&\a(m_{(-1)})\c\mu(l)\o\a(n_{(-1)})\c m_{0}\o\mu(n_0).
\end{eqnarray*}
So we have $R_{M}^{12}\circ R_{M}^{23}=R_{M}^{23}\circ R_{M}^{12}$, as desired.
And this finishes the proof.

\begin{center}
 {\bf ACKNOWLEDGEMENT}
 \end{center}

The work of S. Wang is  supported by  the Anhui Provincial Natural Science Foundation (No. 1908085MA03).
The work of X. Zhang is  supported by the NSF of China (No. 11801304) and the Young Talents Invitation Program of Shandong Province.
The work of S. Guo is  supported by  the NSF of China (No. 11761017)
    and Guizhou Provincial  Science and Technology  Foundation (No. [2020]1Y005).

\renewcommand{\refname}{REFERENCES}

\end{document}